\documentclass[twocolumn, 11pt]{article}
\usepackage{times}
\usepackage[centertags]{amsmath}
\usepackage{amssymb}
\usepackage{amsthm}
\usepackage{amsfonts,amssymb}

\usepackage{mathrsfs}
\usepackage[nodots]{numcompress}
\usepackage{graphicx}

\newtheorem{definition}{Definition}

\newtheorem{lemma}[definition]{Lemma}
\newtheorem{theorem}[definition]{Theorem}
\newtheorem{proposition}[definition]{Proposition}

\DeclareMathOperator{\pd}{\partial}

\DeclareFontFamily{OT1}{pzc}{}
\DeclareFontShape{OT1}{pzc}{m}{it}{<-> s * [1.10] pzcmi7t}{}
\DeclareMathAlphabet{\mathpzc}{OT1}{pzc}{m}{it}

\DeclareMathSymbol{\R}{\mathalpha}{AMSb}{"52}
\DeclareMathSymbol{\C}{\mathalpha}{AMSb}{"43}
\newcommand{\mbb}[1]{\mathbb{#1}}

\newcommand{\K}{\mbb{K}}

\newcommand*{\Scale}[2][4]{\scalebox{#1}{$#2$}}%

\newcommand{\comment}[1]{}

\newcommand{\pzcm}{\mathpzc{m}}
\newcommand{\fq}{\mathfrak{q}}

\newcommand{\bv}{\mathbf{v}}
\newcommand{\bt}{\mathbf{t}}
\newcommand{\bb}{\mathbf{b}}

\newcommand{\bs}{\mathbf{s}}
\newcommand{\bx}{\mathbf{x}}
\newcommand{\by}{\mathbf{y}}
\newcommand{\bz}{\mathbf{z}}
\newcommand{\bw}{\mathbf{w}}

\newcommand{\msC}{\mathscr{C}}

\newcommand{\mcm}{\mathcal{M}}

\newcommand{\mco}{\mathcal{O}}
\begin{document}
\global\def\refname{{\normalsize \it References:}}
\baselineskip 12.5pt
%
%
% TITLE, AUTHOR, ABSTRACT, KEYWORDS
%
\title{\LARGE \bf Generalization of non-Cartan Symmetries to arbitrary dimensions}

\date{}

\author{\hspace*{-10pt}
\begin{minipage}[t]{2.7in} \normalsize \baselineskip 12.5pt
\centerline{J.C. Ndogmo}
\centerline{University of Venda}
\centerline{Department of Mathematics and Applied Mathematics}
\centerline{P/B X5050, Thohoyandou 0950}
\centerline{South Africa}
\centerline{jean-claude.ndogmo@univen.ac.za}
\end{minipage} \kern 0in
\comment{
\begin{minipage}[t]{2.7in} \normalsize \baselineskip 12.5pt
\centerline{SECOND AUTHOR}
\centerline{Name of the University}
\centerline{Institute of Mechanical Engineering}
\centerline{47 West Lincoln Avenue, 87 115 City}
\centerline{COUNTRY}
\centerline{second.author@math.univ.ab}
\end{minipage}
}
%
% If you are three authors then you can use three mini--pages
% instead of two. Their horizontal size must be less than 2.7in
% indicated above. It can be e.g. 2.3in. However, you must pay
% attention that you do not exceed the total width of the text.
%
\\ \\ \hspace*{-10pt}
\begin{minipage}[b]{6.9in} \normalsize
\baselineskip 12.5pt {\it Abstract:}
% The text of the abstract follows.
Second order scalar ordinary differential equations ({\sc ode}s) which are linearizable possess special types of symmetries. These are the only symmetries which are non fiber-preserving in the linearized form of the equation, and they are called non-Cartan symmetries and  known only for scalar {\sc ode}s. We give explicit expressions of non-Cartan symmetries for systems of {\sc ode}s of arbitrary dimensions and show that they form an abelian Lie algebra. It is however shown that the natural extension of these non-Cartan symmetries to arbitrary dimensions is applicable only to the natural extension of scalar second order equations to higher dimensions, that is, to  equivalence classes under point transformations of the trivial vector equation. More precisely, it is shown that non-Cartan symmetries characterize linear systems of {\sc ode}s reducible by point transformation to their trivial counterpart, and we verify  that they do not characterize nonlinear systems of {\sc ode}s  having this property. It is also shown amongst others that the non-Cartan property of a symmetry vector is coordinate-free. Some examples of application of these results are discussed.
\\ [4mm] {\it Key--Words:}
% The key-words follow.
Lie point symmetry algebras, non-Cartan symmetries, Systems of ordinary differential equations
\end{minipage}
\vspace{-10pt}}

\maketitle

\thispagestyle{empty} \pagestyle{empty}
% numbers of pages are supplemented by the editor
%
% THE BEGINNING OF THE TEXT
%
\section{Introduction}
\label{s:intro} \vspace{-4pt}

%%%%%%%%%%%%%%%%%%%%%%%%
%%%%%%%%%%%%%%%%%%%%%%%%
Systems of linear or nonlinear ordinary differential equations ({\sc ode}s) frequently occur in dynamical systems and in many other mathematically based fields, and symmetries are valuable tools for studying such systems. The term symmetry here refers to a generator of the Lie point symmetry algebra of the system of equations. These symmetries provide amongst others valuable means for their identification, the determination of their first integrals and solutions, as well as the qualitative study of these solutions, etc. Many nonlinear systems of {\sc ode}s that occur in real world applications and of course also in theoretical context are actually linear systems in a disguised form, in the sense that such nonlinear systems can be reduced to linear ones by point transformations. Their study is therefore essentially the same as the study of their linear counterpart. We shall consequently focus our attention in this note on linear systems of {\sc ode}s.\par

It is well known that two of the eight Lie point symmetries generating the symmetry algebra of a second order linear {\sc ode} are non fiber-preserving and called non-Cartan. This concept has been known only for scalar equations \cite{ndog17cns, moyoL, chavchanidze}, and we obtain in this paper their natural extension to systems of second order {\sc ode}s of arbitrary dimensions. However, it appears that such extension holds only in equivalence classes (under point transformations) of trivial equations $\by^{(2)}=0,\; \by \in \R^{\pzcm}$ and such classes will be referred to as canonical classes.  This confirms the known fact that the natural extension of scalar linear  {\sc ode}s to systems of arbitrary dimensions consists only of the canonical class. All non-Cartan symmetries are explicitly determined in this paper for arbitrary systems of $\pzcm$ linear equations in canonical classes and it is shown that they form an abelian Lie algebra of dimension $2 \pzcm.$ \par

Moreover, it is shown that non-Cartan symmetries characterize linear systems of {\sc ode}s reducible by point transformation to their trivial counterpart, and we verify  that they do not characterize nonlinear systems of {\sc ode}s  having this property. It is also shown amongst others that the non-Cartan property of a symmetry vector is coordinate-free. Some examples of application of these results are discussed.

%%%%
%%%%

\section{Scalar equations}
\label{s:scalar}\vspace{-4pt}
%%%%%%%%%%%%%%%%%%%%%
%%%%%%%%%%%%%%%%%%%%%

In order to fix ideas, let us consider the Lie point symmetry algebra of the (trivial) free fall equation
\begin{equation} \label{free2}
y''(x)=0,
\end{equation}

 using the  notation $y'= \frac{d y}{d x} = y^{(1)},$ and $y''= \frac{d^2 y}{d x^2} = y^{(2)},$ etc. It is the symmetry algebra of the maximal possible dimension for any second order {\sc ode}, and it would be relevant  to list here its eight generators, given by

%%%
\begin{subequations} \label{symfree2}
\begin{alignat}{3}
S_1 &= \pd_y ,\quad&     S_2 &= x \pd_y,\quad&  F_z &= 2 x \pd_x + y \pd_y,  \\
  F_m &= \pd_x,\quad&  F_p &= x^2 \pd_x + x y \pd_y,\quad&    & \\
H &= y \pd_y,\quad&   C_1  &= y \pd_x,\quad&   C_2 &=  x y \pd_x + y^2 \pd_y.
\end{alignat}
\end{subequations}

The only two non fiber-preserving symmetries in \eqref{symfree2}, i.e. those whose first components do not involve the dependent variable $y,$ are the last two ones, $C_1$ and $C_2,$ and they are thus the non-Cartan symmetries of \eqref{free2}. Moreover, we formally prove in the next proposition that every second order linear {\sc ode} has precisely two non-Cartan symmetries.  First we note that any such equation can always be obtained  from the trivial equation \eqref{free2} through the equivalence transformation
\begin{equation}\label{eqvnh}
x = \rho(t),\quad y= \pi (t) u + \sigma (t),
\end{equation}
which is in fact the most general point transformation leaving invariant any nonhomogeneous linear scalar equation of any given order.

\begin{proposition}\label{p:scalar2}
Every scalar linear second order {\sc ode} has precisely two non-Cartan symmetries.
\end{proposition}

\begin{proof}
Let $\bv =  \xi (x,y) \pd_x + \phi(x,y) \pd_y$  be a symmetry generator of \eqref{free2} and suppose that under a change of variables of the form \eqref{eqvnh} it has expression $\bv = \eta(t,u) \pd_t + \psi(t,u) \pd_u.$ Then $\eta= \eta(t,u)$ clearly satisfies $\eta = \xi \frac{\pd t}{\pd x} + \phi \frac{\pd t}{ \pd y}.$ Since \eqref{eqvnh} is invertible, one must have $t= \mu (x)$ for a certain function $\mu$ depending on $x$ alone, and thus $0 \neq \frac{\pd t}{\pd x}$ is expressible in terms of $t$ alone. Hence it follows from \eqref{symfree2} and \eqref{eqvnh} that  $\eta = \xi (x,y) \frac{\pd t}{\pd x}$ is the first component of a non-Cartan symmetry, that is, it depends explicitly on $u$  if and only if  $\xi$ depends explicitly on $y.$  Therefore, since there are only two non-cartan symmetries in \eqref{symfree2}, the transformed version of \eqref{free2} must also have exactly two non-Cartan symmetries, due to the invertibility of the equivalence transformation \eqref{eqvnh}. This completes the proof of the proposition.
\end{proof}

Let $\Omega= r \frac{d}{d x} + s$ be an ordinary differential operator, where $r$ and $s$ are given functions of $x.$ It is well-known \cite{KM} that a scalar linear homogeneous {\sc ode} of arbitrary  order is reducible by a point transformation to the trivial equation $y^{(n)}=0$ if and only if it is iterative, that is, of the form $\Omega^n[y]=0.$ The normal form of these iterative equations is more convenient for their study, in particular because in such a form they depend on a single arbitrary function. Also, only this form of reduction is always possible in practice without invoking any solution of the equation, and in fact these solutions are often not available even for second order equations.   Here, normal form refers to the form of the {\sc ode} in which the coefficient of the term of second highest order has vanished. Let
\begin{equation} \label{srce}
y'' + \fq y =0
\end{equation}
be the normal form of the second order equation $\Omega^2[y]=0$  for some given and fixed values of $r$ and $s.$ The requirements for the equation to be in normal form forces $s$ to become expressible in terms of $r$ and its derivatives. On the other hand, let

\begin{equation} \label{noreq1}
y^{(n)} + A_n^2\, y^{(n-2)} + \dots + A_n^j\, y^{(n-j)} + \dots + A_n^n\, y=0
\end{equation}
be the normal form of $\Omega^n[y]=0$  corresponding to the same initial parameter  $r$ of $\Omega.$ Then the coefficient $A_n^j$ in \eqref{noreq1}
 are differential polynomials in $\fq$ \cite{JF, ndog17aaa}. Moreover if we let $u$ and $v$ be two linearly independent solutions of \eqref{srce}, then their Wronskian $u v' - u' v$ is a constant which we shall normalize to one. It can then be shown that     $n$ linearly independent solutions of \eqref{noreq1} are given by
\begin{equation} \label{sk}
s_k= u^{n-1-k} v^k, \qquad \text{for $k=0, \dots, n-1.$}
\end{equation}
Every point transformation of the form
\begin{equation} \label{trfsrce}
y= \frac{1}{\lambda} u^{n-1} w,\qquad z= \frac{v}{u}.
\end{equation}
where $\lambda \neq 0$ is a constant reduces any $n$th order equation \eqref{noreq1} to the trivial counterpart $y^{(n)}=0.$ In particular, using \eqref{trfsrce} and \eqref{symfree2} one readily sees that the non-Cartan symmetries corresponding to \eqref{srce} are given by
\begin{equation} \label{nonC22}
C_{11} = yu \pd_x +\, y^2 u' \pd_y,\quad  C_{12} = y v \pd_x +\, y^2 v' \pd_y.
\end{equation}

\section{Generalization to systems of linear equations}
\label{s:transfo} \vspace{-4pt}
%%%%%%%%%%%%%%%%%%%%%%%%%%%%%%%%%%%%%%%%%%%%%%%%
%%%%%%%%%%%%%%%%%%%%%%%%%%%%%%%%%%%%%%%%%%%%%%%%
For an arbitrary   system of {\sc ode}s in $\pzcm$ dependent variables $(y_1, \dots, y_\pzcm)= \by$ and one independent variable $x,$ a symmetry generator
\[
\bv = \xi(x, \by) \pd_x+ \sum_{j=1}^\pzcm \phi_j(x, \by) \pd_j
\]
of the system will be called  non-Cartan if the function $\xi= \xi(x, \by)$ depends explicitly on at least one of the $\pzcm$ dependent variables $y_j.$  Let us denote by $\msC_{\pzcm, n}$ the canonical class of systems of linear {\sc ode}s of order $n$ and dimension $\pzcm.$    The complete algebraic structure of the Lie point symmetry algebra $L_{\pzcm, n}$  of  members of $\msC_{\pzcm, n}$  considered in their most general and normal forms  has recently  been obtained in \cite{ndog17cns} for arbitrary values of $\pzcm$ and for $n \geq 3.$  In fact as a corollary of this result it was shown in the same paper that the maximal dimension of the symmetry algebra for any system of $\pzcm$ (linear  or nonlinear) {\sc ode}s of order $n \geq 3$ is $\pzcm + n \pzcm + 3,$ and that this maximum is achieved precisely  on $\msC_{\pzcm, n}.$  Note that according to a result of \cite{char-nj}, the normal form of a linear systems in $\msC_{\pzcm, n}$ consists of isotropic systems of the form
\begin{equation} \label{norCmn}
\by^{(n)} + A_n^2\, \by^{(n-2)} + \dots + A_n^j\, \by^{(n-j)} + \dots +  A_n^n\, \by=0,
\end{equation}
where $\by = (y_1,\dots, y_\pzcm) \in \R^{\pzcm},$ and the $A_n^j = A_n^j (x)$ are the same scalars appearing in \eqref{noreq1}. In other words, \eqref{norCmn} consists of copies of the same iterative equation. Set $\pd_j= \pd_{y_j}.$ Then a basis of  generators of $L_{\pzcm, n}$ as constructed one by one in \cite{ndog17cns} are given by
\begin{subequations}\label{symLmn}
\begin{align}
H_{i j} &=   y_i \pd_j,\quad \text{ for $i,j = 1, \dots, \pzcm$} \\
S_{k j} & =  s_k \pd_j,\quad \text{for $k= 1,\dots, n$ and $j= 1, \dots, \pzcm $},\\
F_p &=  v^2 \pd_x + (n-1) v v'\, \Scale[0.90]{ \sum_{i=1}^\pzcm} y_i \pd_i,\\
F_m &=  -u^2 \pd_x - (n-1) uu'\, \Scale[0.90]{ \sum_{i=1}^\pzcm} y_i \pd_i \\
F_z &=  2 u v \pd_x + (n-1) (u v' + u' v )\, \Scale[0.90]{ \sum_{i=1}^\pzcm} y_i \pd_i
\end{align}
\end{subequations}
where the $s_k$ are as given by \eqref{sk}. Although the generators \eqref{symLmn} are originally constructed only for Lie algebras $L_{\pzcm, n}$ with $n \geq 3,$ it turns out that by letting $n=2$ in \eqref{symLmn}, the resulting generators are also linearly independent symmetries of $L_{\pzcm, 2}.$ However, by a result of \cite{lopez-slde}, the dimension $d_{\pzcm, 2}$ of $L_{\pzcm, 2}$ is $\pzcm^2 + 4 \pzcm + 3,$ and this can be written as
%%%%
%%%%
\[
d_{\pzcm, 2} = (\pzcm^2 + n \pzcm + 3) + n \pzcm,\quad \text {for $n=2.$}
\]
This shows that in addition to the $\pzcm^2 + 2 \pzcm + 3$ generators of the form \eqref{symLmn}, $L_{\pzcm, 2}$ has precisely $2 \pzcm$ additional generators. The generators of $L_{\pzcm, 2}$ were also obtained in \cite{lopez-slde} but not in a form that exhibit the non-Cartan ones. In fact, there is no reference of any kind to non-Cartan symmetries in \cite{lopez-slde}. It is however clear that the missing $2 \pzcm$ symmetries include all the non-Cartan ones which can be found by calculation for systems of low dimensions. In fact, it turns out that one of the easiest ways to find the $2 \pzcm$ symmetries is to try to guess their expressions from that for scalar equations given by  \eqref{nonC22}. One thus obtain the following result.
\begin{theorem} \label{t:main}
For all $\pzcm \geq 1$ the Lie point symmetry algebra $L_{\pzcm, 2}$ of systems of $\pzcm$ linear {\sc ode}s of order 2 in $\msC_{\pzcm, 2}$ taken in the normal form \eqref{norCmn} has $2 \pzcm$ non-Cartan symmetries
\begin{equation} \label{ncartansym}
\begin{split}
C_{ik} &= y_i u_k \pd_x + \sum_{j=1}^\pzcm y_i y_j u_k' \pd_j, \\
&\text{for  $i= 1, \dots, \pzcm$ and $k=1,2,$}
\end{split}
%%%
\end{equation}
where  $u_k' = d u_k/ dx,$ $u_1=u, u_2 = v,$  and as usual $u$ and $v$ are the two linearly independent solutions of \eqref{srce}. Moreover, these non-Cartan symmetries form an abelian Lie algebra.
\end{theorem}
\begin{proof}
Suppose that $\Delta \equiv (\Delta_1, \dots, \Delta_\pzcm)=0$ is a system of  differential equations defined in the space $\mcm$ of independent and dependent variables, and let $\bv$ be a vector field on $\mcm.$ Then $\bv$ is a symmetry of $\Delta$ if and only if it satisfies the infinitesimal condition of invariance given by
\begin{equation}\label{invcond}
\bv^{(p)} (\Delta_\nu)\big \vert_{(\Delta=0)} =0, \text{ for all $\nu=1, \dots, \pzcm.$}
\end{equation}
where $\bv^{(p)}$ is the $p$th prolongation of $\bv$ to the $p$th jet space of $\mcm.$ Thus the $C_{ik}$ in the theorem are symmetries of $L_{\pzcm, 2}$ because they satisfy the above infinitesimal invariance condition applied to the corresponding second order system of the form \eqref{norCmn}. It is also straightforward to verify that the $C_{ik}$ are pairwise commutative, and thus they form an abelian Lie algebra.
\end{proof}

\section{Characterization of systems of {\sc ode}s admitting non-Cartan symmetries}
\vspace{-4pt}

A question that naturally arises at this point is whether non-Cartan symmetries also exist for  systems of linear {\sc ode}s not belonging to $\msC_{\pzcm, 2},$ that is, which are not members of a canonical class. More generally, the results we have obtained up to now in this paper point to the question of whether non linearizable systems of second order {\sc ode}s may possess non-Cartan symmetries. An answer to the latter question is provided by the trivial system with $\pzcm = 1,$ that is, by scalar {\sc ode}s. It turns out indeed, that some non linearizable second order {\sc ode}s do admit non-Cartan symmetries.  To verify this fact it will be enough to find the general form of scalar second order odes.
\begin{equation}\label{gnl2ode}
y''= F(x, y, y')
\end{equation}
%%%%
admitting the non-Cartan symmetries $C_1= y \pd_x$ and $C_2 = y x \pd_x + y^2 \pd_y$ appearing in \eqref{symfree2}. We first note that \eqref{gnl2ode} admits the symmetry $C_1$ if and only if the function $F$ satisfies
%%%%
\begin{equation}
y F_x - 3 p F + p^2 F_p =0,
\end{equation}
where $p= y'.$ Solving the latter partial differentia equation for $F$ yields
%%%%
\begin{equation}\label{Fexpr}
F= (p/y)^3 G(y, \frac{x}{y} - \frac{1}{p}),
\end{equation}
where $G$ is an arbitrary function of two arguments. With the new expression for $F$ given by \eqref{Fexpr}, it follows that $C_2$ is a also a symmetry of \eqref{gnl2ode} if and only if the function $G$ in \eqref{Fexpr} is of the form
%%%
%%%
\begin{equation*}
G(x,u) = H(x u),
\end{equation*}
for some arbitrary function $H$ of a single argument. Consequently, the most general scalar second order {\sc ode}  admitting the symmetries $C_1$ and $C_2$ is of the form
%%%
%%%
\begin{equation} \label{gnlc1c2}
y''= (p/y)^3 H(x- \frac{y}{p}).
\end{equation}
Consider for instance the particular case of an equation of the form \eqref{gnlc1c2}, given by
\begin{equation}\label{nonlinCart}
y''=   \frac{ p^3 [p (x+1)-y]}{y^3(y- xp)}.
\end{equation}
It clearly follows from Lie's linearization algorithm for second order {\sc ode}s \cite{lielin,meleshlin3rd, invsix} that the latter equation is not linearizable for the obvious reason that is is not a polynomial of degree at most 3 in $p.$ Yet, by construction it admits the symmetries $C_1$ and $C_2.$ We have thus established the fact that non linearizable equations may admit non-Cartan symmetries. It should also be noted that in order to establish this fact, it was enough to exhibit only a non linearizable {\sc ode}  that does admit any given non-Cartan symmetry. In particular, it was enough to consider in the discussion leading to a counterexample of the form \eqref{nonlinCart} only one of the two symmetries $C_1$ and $C_2.$  \par

From the results obtained up to this point in this section, we can affirm that non-Cartan symmetries do not characterize general systems of {\sc ode}s reducible by point transformations to their trivial counterpart, given that some non linearizable systems of {\sc ode}s do admit non-Cartan symmetries. However, it turns out that the non-Cartan symmetries do characterize linear systems of second order {\sc ode}s reducible by point transformation to their trivial counterpart, and this is the case at least for systems of $\pzcm=2$  equations which we shall prove. First, we note that every system of two second order {\sc ode}s can always be put into the normal form

\begin{align} \label{nor2ode}
\by'' &= \begin{pmatrix}  a_1 & a_2\\ a_3 & a_4  \end{pmatrix} \by,\qquad  \by= (y_1, y_2) \in \R^2
\end{align}
where the entries $a_1, a_2, a_3$ and $a_4,$ are arbitrary functions of the independent variable $x.$ By a result of \cite{char-nj}, the equivalence group of \eqref{nor2ode} consists of invertible point transformations of the form
\begin{equation}\label{eqvnor}
x= f(z), \qquad \by = f'(z)^{1/2}\, C\, \bw,
\end{equation}
where $f$ is a  smooth function and $C=(C_{ij}) \in \K^{4}$ is a constant matrix. Applying \eqref{eqvnor} to \eqref{nor2ode}  transforms the latter to an equivalent equation of the form
\begin{align} \label{nor2odev2}
\by'' &= \begin{pmatrix}  A & B\\ C & -A  \end{pmatrix} \by,\qquad  \by= (y, w) \in \R^2
\end{align}
depending on only three arbitrary functions $A, B,$ and $C,$ provided that the function $f=f(z)$  satisfies the nonlinear ode
\begin{equation}\label{eq4f}
-2 [(a_1 + a_4)\circ f]f'^4 - 3 f''^{\,2} + 2 f'f'''=0.
\end{equation}
Let $z= g(x)$ be the inverse of the function $x= f(z).$ By virtue of the invertibility of \eqref{eqvnor}, the function $g$ exists and it follows from \eqref{eqvnor} that the equation satisfied by the auxiliary function $q= g'(x)$ is given by
\begin{equation}\label{eq4q}
-2 (a_1 + a_4)\,q^2 + 3 q'^{\, 2} - 2 q q''=0.
\end{equation}
It follows from Lie's linearization algorithm for second order {\sc ode}s \cite{lielin} that \eqref{eq4q} is linearizable, and hence integrable. Consequently, \eqref{eq4f} is also integrable. Hence without loss of generality, we may assume that any system of two linear second order {\sc ode}s is of the form \eqref{nor2odev2}.\par

Let an  $n$th order system of $\pzcm$ linear equations be given in the form
\begin{align}
\by^{(n)} &+ A_{n-1} \,\by^{(n-1)} + \dots + A_1\,\by ' + A_0 \,\by =\bb, \label{glin} \\
 \by& =(y_1,\dots, y_\pzcm) \in \R^\pzcm,
\end{align}
where the $\pzcm \times \pzcm$ matrices  $A_j= A_j(x)$ are given functions of the independent variable $x,$ and $\bb = (b_1(x), \dots, b_\pzcm (x))$ is the nonhomogeneous  term. We recall here that by a result of \cite{char-nj} the equivalence group of \eqref{glin} consists of invertible point transformations of the form
%%%
\begin{equation}\label{eqvlin}
x= f(z), \qquad \by = Q \bw + \bs,\qquad \bw, \bs \in \R^\pzcm,
\end{equation}
where $Q=Q(z) =\big(q_{ij}\big)$ is an $\pzcm \times \pzcm$ matrix and $\bs= \bs(z)= (s_1,\dots, s_\pzcm)$ is a particular solution of \eqref{glin}.
%%%
\begin{lemma}\label{le:chgsym}
The non-Cartan property of a symmetry of a given system of linear {\sc ode}s is coordinate-free. That is, a symmetry generator of a linear system of {\sc ode}s is non-Cartan in a given coordinate system if and only if it remains non-Cartan under the general change of variables \eqref{eqvlin}.
\end{lemma}
\begin{proof}
Let $\mco$ be an open subset of $\R\times\R^\pzcm$ coordinatized by $\bx = (x, \by)$  where $\by=(y_1, \dots, y_\pzcm).$ Let
%%%
\begin{equation} \label{gnlsym}
\bv= \xi(\bx) \pd_x + \sum_i \phi_i(\bx) \pd_j
\end{equation}
be a vector field on $\mco.$  Denote by $\bz= \psi(\bx)$ a change of coordinates of the general form \eqref{eqvlin}. Thus $\bz= (z, \bw),$ with $\bw= (w_1,\dots, w_\pzcm)$.  To prove the lemma it suffices to show that $\bv$ is non-Cartan in the $\bx$ coordinates if and only if it is non-Cartan in the $\bz$ coordinates.  Suppose that in the $\bz$ coordinates $\bv$ has expression
\[
\bv = \bar{\xi}(\bz) \pd_z + \sum_i \bar{\phi}_i(\bz) \pd_{w_j}.
\]

It is well known from standard results on transformation groups \cite{olver-EIS, sagle} that in the $\bz$ coordinates, the components $\bar{\xi}$ and $\bar{\phi_i}$ of $\bv$ are given by
%%%%%%%%
\begin{subequations}\label{chgsym}
\begin{align}
\bar{\xi}(\bz) &= \xi(\bx) \frac{\pd \psi_0}{\pd x} + \sum_{i=1}^\pzcm \phi_i(\bx) \frac{\pd \psi_0}{\pd y_i} \label{chgsym1} \\
\bar{\phi}_j(\bz) &= \xi(\bx) \frac{\pd \psi_j}{\pd x} + \sum_{i=1}^\pzcm \phi_i(\bx) \frac{\pd \psi_j}{\pd y_i},
\end{align}
\end{subequations}
where $(\psi_0, \psi_1,\dots, \psi_\pzcm)=  \psi.$ Since $\bz= \psi(\bx)$ represents the change of coordinates \eqref{eqvlin} one has $z= \psi_0(\bx)= g(x),$  where $g$ is the inverse of the function $f$ from \eqref{eqvlin}. Consequently, it follows from  \eqref{chgsym1} that
%%%%
%%%%
\begin{subequations} \label{bvj's}
\begin{align}
\bar{\xi}(\bz) &= \frac{1 }{f'(z)}  \xi ( \bx  )=   \frac{1 }{f'(z)}  \xi (f(z), y_1, \dots, y_\pzcm) \\[-2mm]
\intertext{where by \eqref{eqvlin} one has \vspace{-3mm}}
y_i &= \sum_{j} q_{ij} w_j + s_i(z).
\end{align}
\end{subequations}
It thus follows from the invertibility of the matrix $Q= \big( q_{ij} \big)$ that each of the variables $y_i$ depends explicitly on at least one of the new dependent variables $w_j.$ Consequently, \eqref{bvj's} clearly shows that $\xi(\bx)$ depends explicitly on one of the original dependent variables $y_i$ if and only if $\bar{\xi}(\bz)$ also depends explicitly on one of  new independent variables $w_j.$ In other words, $\bv$ is non-Cartan in the $\bx$ coordinates if and only if it is non-Cartan in the $\bz$ coordinates and this completes the proof of the lemma.
\end{proof}

It should  be noted that Proposition \ref{p:scalar2} may also be obtained as a corollary of Lemma \ref{le:chgsym}.

\begin{theorem} \label{t:id-ncart}
An arbitrary  linear system (S) of two second order {\sc ode}s is reducible by a point transformation to the trivial equation $\by''=0, \by \in \R^2$ if and only if it admits a non-Cartan symmetry.
\end{theorem}

\begin{proof}
It follows from Lemma \ref{le:chgsym} that without loss of generality one may assume that the given system of {\sc ode}s is in the reduced normal form \eqref{nor2odev2}. Moreover, by Theorem \ref{t:main} and Lemma \ref{le:chgsym} it suffices to show that if an equation of the form \eqref{nor2odev2} admits a non-Cartan symmetry, then it is trivial. Therefore, let a vector field of the form \eqref{gnlsym} with $\pzcm=2$ be a symmetry generator of \eqref{nor2odev2}. By applying the second prolongation of $\bv$ to \eqref{nor2odev2} according to the infinitesimal invariance criterion  \eqref{invcond} and then expanding the resulting expression as polynomials in the derivatives of $y(z)$ and $w(z)$ yields the so-called determining equations for $\bv.$ In this instance they are given by

%%%
%%%
\begin{subequations}\label{deq1}
\begin{align}
&\xi _{ww}=0,\quad  \xi _{yw}=0, \quad \xi _{yy}=0   \label{deq1n1}\\
&\eta _{ww}=0,\qquad  \phi _{yy}=0 \label{deq1n2}\\[1.5mm]
&\phi _{ww}-2 \xi _{xw}=0,\qquad 2 \eta _{yw}-2 \xi _{xw}=0  \label{deq1n3}\\
&\eta _{yy}-2 \xi _{xy}=0,\qquad  -2 \xi _{xy}+2 \phi _{yw}=0 \label{deq1n4}\\[1.5mm]
&  2 \eta _{xw}-2 B w \xi _w-2 A y \xi _w=0 \label{deq1n5}\\
& 2 A w \xi _y-2 C y \xi _y+2 \phi _{xy}=0 \label{deq1n6}\\[1.5mm]
\begin{split}
 &-A \eta -B \phi -y \xi  A_x-w \xi  B_x-A w \eta _w+C y \eta _w\\
 &+B w \eta _y+A y \eta _y-2 B w \xi _x-2 A y \xi _x+\eta _{xx}=0
\end{split}\\
\begin{split}
& 3 A w \xi _w-3 C y \xi _w-B w \xi _y\\
&-A y \xi _y-\xi _{xx}+2 \phi _{xw}=0
\end{split}\\
\begin{split}
& A w \xi _w-C y \xi _w-3 B w \xi _y\\
&-3 A y \xi _y+2 \eta _{xy}-\xi _{xx}=0
\end{split}\\
 \begin{split}
&-C \eta +A \phi +w \xi  A_x-y \xi  C_x+2 A w \xi _x-2 C y \xi _x \\
&-A w \phi _w+C y \phi _w+B w \phi _y+A y \phi _y+\phi _{xx}=0.
 \end{split}
\end{align}
\end{subequations}
%%%
It thus follows from \eqref{deq1n1} that
\begin{equation} \label{xin1}
\xi = \alpha y + \beta w + \gamma,
\end{equation}
where $\alpha,$  $\beta,$ and $\gamma$ are some functions of $x.$ It also follows from \eqref{deq1n2} that
\[
\eta= a w + b, \qquad \text{ and }\qquad  \phi = R y + S,
\]
where $a$ and $b$ are functions of $x$ and $y$ while $R$ and $S$ are functions of $x$ and $w.$ Substituting these new expressions for $\xi, \eta,$  and $\phi$ into \eqref{deq1} and solving the resulting version of \eqref{deq1n3} and \eqref{deq1n4} shows that
%%%
\begin{subequations} \label{newRS}
\begin{align}
R&= \alpha' w + r_1,\qquad S= \beta' w^2 + s_1 w + s_2\\
a&= \beta' y + a_1,\qquad b= \alpha' y^2 + b_1 y + b_2,
\end{align}
\end{subequations}
where $r_1, s_1, s_2, a_1, b_1,$ and $b_2$ are of course some functions of $x.$ Substituting the new expressions from \eqref{newRS} into the latest version of \eqref{deq1}, it follows from the resulting version of \eqref{deq1n5} and  \eqref{deq1n6} that $r_1= k_1$ and $a_1= k_2$ for some constants $k_1$ and $k_2.$ Updating the latest version of \eqref{deq1} with these new values for $r_1$ and $a_1$ yields the new system of remaining determining equations

%%%
\begin{subequations} \label{deq2}
\begin{align}
&-2 C y \alpha +2 w \left(A \alpha +\alpha _{xx}\right)=0  \label{deq2n1} \\
 &-2 B w \beta -2 y \left(A \beta -\beta _{xx}\right)=0  \label{deq2n2}\\[2mm]
 \begin{split}
& -A b_2-B S_2+w^2 \left(-\beta  B_x-2 B \beta _x\right)\\
&+w \left(-2 A k_2+B b_1-\gamma  B_x-B S_1-2 B \gamma _x\right)\\
&+{b_2}_{x,x} +y^2\left(-\alpha  A_x-A \alpha _x+C \beta _x+\alpha _{xxx}\right)\\
&+y \big(-B k_1+k_2 C-\gamma  A_x-2 A \gamma _x+  {b_1}_{xx}\\
&+w \left(-\beta  A_x-\alpha  B_x-B \alpha _x-3 A \beta _x+\beta _{xxx}\right)\big)=0 \label{deq2n3}
\end{split}\\[2mm]
\begin{split}
&2 {s_1}_x+y \left(-A \alpha -3 C \beta +\alpha _{xx}\right)\\
&+w \left(-B \alpha +3 A \beta +3 \beta _{xx}\right)-\gamma _{xx}=0 \label{deq2n4}
\end{split}\\[1mm]
\begin{split}
&2 {b_1}_x+y \left(-3 A \alpha -C \beta +3 \alpha _{xx}\right)\\
&+w \left(-3 B \alpha +A \beta +\beta _{xx}\right)-\gamma _{xx}=0 \label{deq2n5}
\end{split}\\[1mm]
\begin{split}
& -C b_2+A S_2+y^2 \left(-\alpha  C_x-2 C \alpha _x\right)\\
&+w \left(B k_1-k_2 C+\gamma  A_x+2 A \gamma _x+{s_1}_{xx}\right) \\
&+{S_2}_{xx}+y \big(2 A k_1-C b_1-\gamma  C_x+C S_1-2 C \gamma _x \\
&+w \left(\alpha  A_x-\beta  C_x+3 A \alpha _x-C \beta _x+\alpha _{xxx}\right)\big)\\
&+w^2 \left(\beta  A_x+B \alpha _x+A \beta _x+\beta _{xxx}\right)=0.
\end{split}
\end{align}
\end{subequations}
Recall that from the expression for $\xi$ given in \eqref{xin1}, $\bv$ is non-Cartan if and only if $\alpha \neq0$ or $\beta \neq 0.$ \par
We first suppose that $\alpha \neq 0.$ Then, by the vanishing of the coefficient of $y$ in \eqref{deq2n1}, it follows that $C=0.$ Substituting this new value for $C$ in \eqref{deq2} and comparing the coefficient of $w$ in the resulting version of \eqref{deq2n4} and \eqref{deq2n5} shows that $B=0$ must hold. With this new value of $B,$ comparing the coefficient of $w$ in \eqref{deq2n1} and the coefficient of $y$ in \eqref{deq2n4} shows that $A=0$ must also hold. In view of the form of the original equation \eqref{nor2odev2}, we have thus shown that if $\alpha \neq0$ then \eqref{nor2odev2} reduces to the trivial equation. \par
We now consider the case $\beta \neq 0.$ Since the roles of $\alpha$ and $\beta$ in the expression of $\xi$ in \eqref{xin1} are clearly symmetrical, it also follows that $\beta \neq 0$ implies that \eqref{nor2odev2} is trivial. In other words, if the symmetry vector $\bv$ is non-Cartan, then the corresponding equation \eqref{nor2odev2} is  trivial, and this completes the proof of the theorem.
\end{proof}

It should be noted that it is not essential in the proof of Theorem \ref{t:id-ncart} to assume that the given system of {\sc ode}s is of the form \eqref{nor2odev2}. One could as well assume even more simply that it is rather of the more general form \eqref{nor2ode}. The proof of the theorem is then carried out in that case along the same lines as in the given proof, except that the resulting system is not  a trivial one, but rather an isotropic one. One can then make use of a result of \cite{char-nj} to conclude that such an isotropic system is reducible to the trivial equation by a point transformation.

\section*{Concluding remarks}
\vspace{-4pt}

We have proved in this paper that the non-Cartan property of a symmetry  is coordinate-free for a system linear {\sc ode}s of arbitrary order and dimension, and that non-Cartan symmetries do not embed a characterization of nonlinear systems reducible to their trivial counterpart by a point transformation. We have however shown that non-Cartan symmetries do characterize linear systems of second order {\sc ode}s  reducible by a point transformation to the trivial system. More exactly, we have provided the proof to this fact for systems of two equations. Although it is apparent that the result holds for any system of $\pzcm \geq 2$ equations, it remains an open problem to systematically prove this fact.\par

An immediate application of this result is an easy method for identifying linear systems of {\sc ode}s reducible by point transformation to their trivial counterpart. Indeed, the results of the paper suggest that one simple way to determine if a given system of linear equations is in $\msC_{\pzcm, n}$ is to find out if it has any non-Cartan symmetry, in case the Lie point symmetry algebra is available. For instance if we let $A= \left(  \begin{smallmatrix}1& 0\\2 & 1 \end{smallmatrix}\right),$ and $\by= (y_1, y_2),$ then the system of linear equations $\by'' + A \by =0$ has a seven dimensional Lie point symmetry algebra, all of whose generators are fiber-preserving. Therefore the given system is not a member of $\msC_{\pzcm, 2}.$\par

When the Lie point symmetry algebra is not available, the other simple way to find out if a given system of linear equations is a member of $\msC_{\pzcm, n}$ is to test if the system admits any  non-Cartan symmetry, and in so doing the component $\xi$ in \eqref{gnlsym} should  be taken in the form \eqref{xin1}.

\end{document}